\documentclass{mpri}
\usepackage{authblk}

\usepackage[ngerman,english]{babel}
\usepackage[latin1]{inputenc}
\usepackage{csquotes}
\usepackage{amsfonts,amsmath}
\usepackage{empheq}
\usepackage[titletoc,title]{appendix}
\usepackage[backend=bibtex,url=false,doi=false,eprint=false,giveninits=true,isbn=false,style=numeric-comp,maxnames=99]{biblatex}
\makeatletter
\def\blx@maxline{77}
\makeatother

\bibliography{bibl_review.bib}
\AtEveryBibitem{\clearlist{language}}

\usepackage{cases}
\usepackage{mathabx}
\usepackage{bbm}
\usepackage{xfrac}
\usepackage{fancyhdr}
\usepackage{color}
\usepackage[colorinlistoftodos,textsize=small,backgroundcolor=white,bordercolor=blue,linecolor=blue,disable]{todonotes}
\usepackage[colorlinks]{hyperref}
\definecolor{blue75}{rgb}{0,0,.75}
\definecolor{green75}{rgb}{0,.75,0}
\hypersetup{colorlinks=true, urlcolor=blue75,linkcolor=blue75,citecolor=green75,pdfstartview=FitB,bookmarksopen=true,bookmarksopenlevel=1}
\usepackage{constants}
\newcommand{\parenthezises}[1]{\arabic{#1}}
\newconstantfamily{C}{
symbol=C,
format=\parenthezises,
reset={section}
}
\newconstantfamily{M}{
symbol=M,
format=\parenthezises,
reset={section}
}
\newconstantfamily{B}{
symbol=B,
format=\parenthezises,
reset={section}
}
\newconstantfamily{xi}{
symbol=\beta,
format=\parenthezises,
reset={section}
}
\usepackage{graphicx}
\graphicspath{{images/} }
\usepackage{wrapfig}
\usepackage{figbib}
\allowdisplaybreaks
\usepackage[capitalise]{cleveref}

\crefdefaultlabelformat{{#2#1#3}}

\crefname{equation}{}{}

\crefname{enumi}{}{}
\creflabelformat{enumi}{{(#2#1#3)}}

\crefname{section}{{ Section}}{{ Sections}}
\crefname{subsection}{{ Subsection}}{{ Subsections}}
\crefname{subsubsection}{{ Paragraph}}{{ Paragraphs}}

\crefname{Theorem}{{ Theorem}}{{ Theorems}}

\newremark{Example}[Theorem]{Example}
\crefname{Example}{{ Example}}{{ Examples}}
\newdefinition{Notation}[Theorem]{Notation}
\crefname{Notation}{{ Notation}}{{ Notations}}

\begin{document}

\newcommand{\cb}{}
 \newcommand{\red}[1]{\textcolor{red}{#1}}
\newcommand{\cmg}[1]{\textcolor{magenta}{#1}}
\newcommand{\cgr}[1]{{#1}}
\newcommand{\D}{\mathbb{D}}
\newcommand{\E}{\mathbb{E}}
\newcommand{\G}{\mathbb{G}}
\newcommand{\T}{\mathbb{T}}
\newcommand{\PP}{\mathbb{P}}
\newcommand{\R}{\mathbb{R}}
\newcommand{\N}{\mathbb{N}}
\newcommand{\F}{\mathbb{F}}
\newcommand{\V}{\mathbb{V}}
\newcommand{\ve}{\varepsilon}
\newcommand{\ce}{{c^{\varepsilon}}}

\def\diam{\operatorname{diam}}
\def\dist{\operatorname{dist}}
\def\diver{\operatorname{div}}
\def\ess{\operatorname{ess}}
\def\inner{\operatorname{int}}
\def\osc{\operatorname{osc}}
\def\sign{\operatorname{sign}}
\def\supp{\operatorname{supp}}
\newcommand{\BMO}{BMO(\Omega)}
\newcommand{\LOne}{L^{1}(\Omega)}
\newcommand{\LOnen}{(L^{1}(\Omega))^d}
\newcommand{\LTwo}{L^{2}(\Omega)}
\newcommand{\Lq}{L^{q}(\Omega)}
\newcommand{\Lp}{L^{2}(\Omega)}
\newcommand{\Lpn}{(L^{2}(\Omega))^d}
\newcommand{\LInf}{L^{\infty}(\Omega)}
\newcommand{\HOneO}{H^{1,0}(\Omega)}
\newcommand{\HTwoO}{H^{2,0}(\Omega)}
\newcommand{\HOne}{H^{1}(\Omega)}
\newcommand{\HTwo}{H^{2}(\Omega)}
\newcommand{\HmOne}{H^{-1}(\Omega)}
\newcommand{\HmTwo}{H^{-2}(\Omega)}

\newcommand{\LlogL}{L\log L(\Omega)}

\def\avint{\mathop{\,\rlap{-}\!\!\int}\nolimits} 

\newcommand{\om}{\omega }
\newcommand{\Om}{\Omega }

\newtheorem{proofpart}{Step}
\makeatletter
\@addtoreset{proofpart}{Theorem}
\makeatother
\numberwithin{equation}{section}
\title[Flux limitation mechanisms in multiscale modelling of cancer invasion]{Flux limitation mechanisms arising in multiscale modelling of cancer invasion}
\author[Zhigun]{\ns A. Zhigun$^*$\par
School of Mathematics and Physics, Queen's University Belfast, University Road, Belfast BT7 1NN, Northern Ireland, UK
}
\email{A.Zhigun@qub.ac.uk}

\titlefootnote{
\\
doi:
 }

\date{9 March 2022}
\maketitle
\begin{abstract}
Tumour invasion is an essential stage of cancer progression. Its main drivers are diffusion and taxis, a directed movement along the gradient of a stimulus. Here we review models with flux limited diffusion and/or taxis which have applications in modelling of cell migration, particularly in cancer. Flux limitation ensures control upon propagation speeds, precluding unnaturally quick spread which is typical for traditional parabolic equations. 
We recall the main properties of models with flux limitation effects and discuss ways to construct them, concentrating on multiscale derivations from kinetic transport equations. 
\end{abstract}
\section{Introduction}\label{sec:intro}
Migration is an essential stage of cancer progression \cite{hanahan}. It starts with  cells of a growing  malignant tumour invading the surrounding tissue matrix. Locally, this eventually leads to malfunctioning of the organ in which the tumour has arisen.  An even more life-threatening implication of invasion is metastasis. Migrating cancer cells are able to reach and  penetrate blood and/or lymph vessels. If this occurs, and the cells manage to survive  the  transportation across the  circulatory system,  they can colonise distant sites,  forming further neoplasms there. This process is termed metastasis. It is responsible for about 90\% of all  deaths coursed by cancer  \cite[chapter 14]{Weinberg}. 

Cancer invasion results 
from a complicated interplay of many effects, including cell movement, proliferation, and interaction among themselves and with their surroundings. Cancer cell motility, proliferation, even survival, are  subject to cell-tissue interaction, see e.g. \cite{pickup}. The main mechanisms  of cancer cell  movement are diffusion and, unlike lifeless particles, various taxes. Taxis refers to movement guided by the gradient of a stimulus in the cell's surroundings. 
One speaks, e.g.  of haptotaxis, chemotaxis, and pH-taxis if the motion is directed by the gradients of tissue fiber density, a diffusing chemical, and pH, respectively. All three kinds of taxis occur during tumour invasion. Haptotaxis plays the key role, directing the cells along the tissue fibers   \cite{carter}. 

Many mathematical models of  invasion have been derived and studied  with the aim of  improving our understanding  of the involved biological phenomena. 
Macroscopic reaction-diffusion-taxis (RDT) systems are among the most popular tools in this context. Yet few of them, even if carefully derived, do not violate a very basic property, namely that  cell speeds are bounded by a certain finite intrinsic value. This value is, for example, independent of the initial cell density. 
Models with flux limitation have been designed with the aim to have control on  speeds as well as other characteristics of propagation. In these models, diffusion and/or taxis parts of the cell flux possess a priori  bounds that are independent of the spacial gradients of the involved quantities. This  generally leads to a finite and well-controlled contribution from the corresponding motion effect (i.e. diffusion or taxis)  to the propagation speed, staying  below a certain value that can be directly determined from the equation coefficients. This property makes equations with flux-limitation an attractive tool for modelling cancer invasion. Some models involving them have been proposed in \cite{conte2021modeling,KumarSur,DKSS20,KimFried}, and we believe that there are more to come.

RDT equations are often obtained by macroscopic flux balance. Examples in the context of cancer invasion include  the models in  \cite{anderson2000,ChapLol2011}. However, this method is not very accurate as it disregards important  information from smaller scales. In contrast, a multiscale approach that is based on construction and upscaling of a mesoscopic kinetic transport equation (KTE)  leads to a considerably more  precise description on the macroscale. An early example of an application to modelling of a type of cell motion often observed in cancer can be found in \cite{HillenM5}.

\bigskip
In this paper we review models with flux-limited (FL) diffusion and/or taxis in the context of cell migration, cancer invasion being the principal application that we have in mind.  
We start with FL diffusion equations in \cref{SecDiff}. We discuss the rational behind such models, their main properties, the analytical challenge that they pose, and how they can be derived on the macroscale. In the same spirit we  overview  RDT systems with FL taxis in \cref{SecTax}.   \cref{Sec:Multi} is the main part of this paper. There we closely examine several ways to derive equations involving FL motion from  KTEs by means of a suitable upscaling.  




\section{FL diffusion equations}\label{SecDiff}
\subsection{Motivation}\label{SecDiffMot}
Early PDE models for population spread mostly include the standard linear diffusion, see e.g.  \cite{MurrayI,MurrayII}. The simplest possible model of this kind is the linear diffusion equation
\begin{align}
 \partial_t c=\nabla_x\cdot(D_c\nabla_x c),\label{LD}
\end{align}
where  $c$ is the density of a population and $D_c>0$ is its diffusion coefficient which is assumed to be constant. However, in many biological applications, the choice of a diffusion flux which is in a constant proportion to the density gradient turns out to be inadequate due to certain well-known characteristics of patterns it induces, including:
\begin{enumerate}
 \item\label{LD1} infinite propagation speed, i.e. even if $c_0$ has a compact support, $c(t,\cdot)$ is positive everywhere  for any $t>0$;
 \item\label{LD2} formation of  smooth Gaussian-like structures, with any initial singularities  eliminated;
 \item\label{LD3} insensitivity to overcrowding, so that in local regions where the population is dense but rather evenly distributed diffusion is not enhanced.
\end{enumerate}
Concerns were  raised  in connection with biofilm formation \cite{eberl2001new}, tumour invasion  \cite{ZSU,conte2021modeling}, as well as more generally     \cite{BBNS2010} (see also  references therein), pointing out that experimentally observed patterns of a spreading cell  population violate  \cref{LD1,LD2,LD3}. 

 A partial remedy is, e.g.   the  diffusion term proposed in \cite{eberl2001new}: 
\begin{align}
 \nabla_x\cdot\left(D_c\frac{c^b}{(c_{max}-c)^{a}}\nabla_x c\right),\label{DegSingD}
\end{align} 
where $a,b, c_{max},D_c>0$ are some intrinsic parameters of the population, $c_{max}$ being the  physically maximal possible density. 
Originating from the porous media equation
\begin{align}
 \partial_tc=\nabla_x\cdot(c^m\nabla_x c)\label{PM}
\end{align}
for a constant $m>0$, the power-like degeneracy at zero enforces  finite speed of propagation. This means that for compactly supported initial data the speed with which the spacial support of the solution extends with time  is finite. 
As to the singularity at $c_{max}$, it serves to  accelerate the  
dispersal from  densely  populated areas where it comes close to that value. 
In this situation one can meaningfully speak of a moving front, i.e. the boundary of the support of $c(t,\cdot)$, and how it changes over time if one starts with an initial datum which is compactly supported. 
Still, using diffusion    \cref{DegSingD} in an   equation for cell motion cannot guarantee an accurate description of a spreading out. Indeed, along  the moving front  solutions   behave  similar to those  of the  porous media equation \cref{PM},  in particular (see e.g.  \cite{Vazquez}):
\begin{enumerate}
 \item the  propagation speed, though finite due to a degeneracy at zero, is  not a population intrinsic  trait. It  
depends on the initial density;
\item along the moving front, the power-like degeneracy  smooths  any sharp singularities, e.g. jump discontinuities   down to  H\"older continuity.
\end{enumerate}
  Thus, with such diffusion as in \cref{DegSingD}, impossibly large speeds are achievable.  Furthermore,  sharp moving fronts cannot be reproduced, and yet this effect has been  observed in experiments, e.g. for glioblastoma spread \cite{conte2021modeling}. We refer to  \cite{CCCSS15} where  further properties of the standard porous media equation are surveyed in relation to modelling of propagation of moving fronts. 

In all mentioned cases  the diffusion flux is of the form 
\begin{align}
 J_{diff}=cV_{diff},\qquad V_{diff}=-b_{diff}(c,\nabla_xc)\label{V}
\end{align}
for some (vector-valued) function $b_{diff}$ which is linear in its second argument and hence unbounded. As discussed  above, for the corresponding diffusion equations there are also no implicit universal bounds, i.e. such that would hold for all solutions, on the propagation speed.  FL diffusion or flux-saturated
diffusion, as it is often called, offers a means for an explicit control on speeds. It corresponds to the situations where for each fixed $c$ the speed function, i.e. $|b|$ 'saturates' in its second argument, meaning that it converges to a finite value as $|y|\to\infty$. An early example is the relativistic heat equation
\begin{align}
 \partial_tc=\nabla_x\cdot\left(D_c\frac{c\nabla_x c}{\sqrt{c^2+\frac{D_c^2}{C^2}|\nabla_x c|^2}}\right)\label{RH}
\end{align}
for constants $D_c,C>0$. 
Here
\begin{align*}
 b_{diff}(c,y)=D_c\frac{y}{\sqrt{c^2+\frac{D_c^2}{C^2}|y|^2}},
\end{align*}
so that for each $c$
\begin{align*}
 |b_{diff}(c,y)|\underset{|y|\to\infty}{\nearrow}C,
\end{align*}
which implies that the saturation condition is satisfied and  that formally cell speeds do not exceed $C$. 
Various modifications of \cref{RH} exist in the literature, including  hybrid models which combine FL diffusion with \cref{PM}  \cite{CCCSS15} or volume saturation effects   \cite{BuriniChouhad}, see also references in these papers.

FL diffusion models were extensively reviewed in \cite{CCCSS15}. There one  may find a historical account on these models and a detailed discussion of their various properties that can be observed in numerical simulations and to a large extent  also verified by rigorous analysis, as well as of the ways they can be derived. 
In particular, it was proved for    \cref{RH} in \cite{ACMM06},
as well as for a broad class of its variants in \cite{Calvo2015} that in contrast to \cref{PM}:
\begin{enumerate}
\item the propagation speed is  bounded by a universal constant which is an explicit model parameter. Moreover, generically this speed is equal to that constant.   For the  relativistic heat equation, the propagation speed is essentially equal to $C$;
\item initial discontinuities on the support boundary are, at least in certain cases, propagated eternally. Thus, in general regularisation does not  occur on the moving front. 
\end{enumerate}
If a reaction term of a Fisher-Kolmogorov type is included into \cref{RH}, then singular fronts can also be propagated, e.g. by  travelling waves, see \cite{CCCSS15} and references therein.
\subsection{Macroscopic derivation}
In this Subsection we briefly review two methods  of deriving FL diffusion on the macroscale. Multiscale alternatives are addressed in \cref{Sec:Multi}.
\subsubsection{Flux adjustment.}
The simplest construction goes back to  \cite{rosenau1992tempered} and consists of  adjusting the form of the diffusion flux directly on the macroscale: one replaces \cref{V} by 
\begin{align*}
  V_{diff}=-\tilde{b}_{diff}(c,\nabla c),\qquad \tilde{b}_{diff}=\psi\circ b_{diff},
\end{align*}
for some (vector-valued) continuous bounded function  $\psi$ that  saturates to a constant $C$ at infinity and is  close to identity for  $|V_{diff}|\ll C$. For example, taking $b_{diff}(c,y)=D_cy$ and 
\begin{align}
\psi(z)=\frac{z}{\sqrt{1+\frac{|z|^2}{C^2}}}\label{psi}
\end{align}
yields \cref{RH}. 
More general modifications of the  formula for $V$ are, of course, possible, e.g. one could take $\psi=\psi(x,z)$ in order to account for local heterogeneity of the surroundings. 
 
In the context of modelling cell migration,  a purely macroscopic framework was used, e.g. in \cite{conte2021modeling} in order to describe the moving fronts observed in  glioblastoma  invasion. 
On the whole, this approach is  flexible, yet  may lead to inaccurate descriptions, see the discussion in \cref{Sec:Multi}.

\subsubsection{Optimal transport.}
An alternative macroscopic derivation can be accomplished with the optimal transport approach. It was noticed in \cite{Brenier03} and made rigorous in \cite{McCP} that in the Monge-Kantorovich mass transportation  framework equation  \cref{RH} is the 
gradient flow of the Boltzmann entropy
\begin{align*}
 F(r)=r\ln r-r
\end{align*}
for the Wasserstein metric  corresponding to the cost function
\begin{align}
 k(z)=\begin{cases}
       C^2\left(1-\sqrt{1-\frac{|z|^2}{C^2}}\right)&\text{if }|z|\leq C,\\
       +\infty&\text{if }|z|>C.
      \end{cases}\label{cost}
\end{align}
 Choosing $k$ differently allows to obtain other variants of FL models, see e.g. examples in \cite{CCCSS15}. For $|z|\ll C$, the cost function in \cref{cost} is close to the quadratic function $\frac{1}{2}|z|^2$. Choosing $k(z)=\frac{1}{2}|z|^2$ for all $z\in\R$ would yield the heat equation. It is the case for which this method was originally  proposed and carried out in   \cite{RKO}. 

\subsection{Analytical challenges}
While the ability of  equations with a FL diffusion  to  reproduce  discontinuous moving  fronts  is attractive for modelling purposes, their potential presence leads to substantial analytical difficulties. Indeed, for \cref{RH} one can only expect that   $u(t,\cdot)$ belongs to the space of functions of bounded variation,  so that its spacial derivatives  are  Radon measures and, in general, not some integrable functions. 
This makes the diffusion flux particularly difficult to handle because it is a nonlinear function of the spacial gradient of $u$. A well-posedness theory for \cref{RH} and its variants was developed and studied in a series of works  \cite{ACMM10,ACM04,ACM05,ACM05Cauchy,ACM08}, as well as  \cite{ACMFK10} which treats a reaction-diffusion equation 
(see also those references in  \cite{CCCSS15} which deal with further modifications of the model).  There a suitable form of a so-called entropy solution was developed and its  existence and uniqueness was proved.

\section{RDT systems with FL mechanisms}\label{SecTax}
 \subsection{Motivation} 
Early PDE models for taxis were developed specifically for chemotaxis. This is a directed movement of cells or organisms in response to diffusing chemical cues. Many  biological processes, including cancer invasion,  
crucially depend on it, see e.g.  \cite{eisenbach2004}.  
PDE systems  modelling chemotaxis
have enjoyed great popularity ever since the introduction of the classical  Keller-Segel model \cite{KS70,KS71}:
\begin{subequations}\label{KS}
 \begin{align}
  &\partial_t c=\nabla_x\cdot(D_c\nabla_x c-\chi c\nabla_x S),\label{KSc}\\
  &\partial_t S=D_v\Delta_x S- \alpha S+ \beta c\label{KSv}
 \end{align}
\end{subequations}
for some constants $\alpha,\beta,D_c,D_v,\chi>0$. 
 Equation \cref{KSc} for the population density $c$  features two motion effects: a linear diffusion and  drift in the direction of the spacial gradient of the concentration of a chemical  $S$. This model is able to reproduce formation of aggregates which is the main implication of chemotaxis. Yet the resulting patterns are often inadequate because one observes:
 \begin{enumerate}
  \item that already in finite time an unlimited aggregation may occur, leading to a so-called 'blow-up' (i.e.  the cell density becomes unbounded); 
  \item consequences of the linear diffusion, see \cref{SecDiff};
  \item an unlimited response to chemotaxis due to the chemotaxis flux being directly proportional to the gradient of the attractant.
 \end{enumerate}
  Arguably the main drawback of \cref{KS} is that it cannot maintain a reasonable balance between 
the two drivers of cell spread. Indeed, there are essentially two options: either the cell motion is  governed by the linear diffusion, and then chemotaxis hardly plays any role, or chemotaxis dominates, inducing an unrealistically strong aggregation, even a blow-up. This aspect is  particularly well-understood,  see e.g. reviews  \cite{Horstmann,BBTW,LankWink2020}. 

Similar to the purely diffusion case, one could try to improve the model by allowing the population diffusion coefficient, $D_c$ and the so-called chemotactic sensitivity, $\chi$ to depend on $c$ and/or $S$, see e.g. \cite{Horstmann,BBTW,HillenPainter} where many examples can be found. 
One such model which includes  diffusion of the form given by \cref{DegSingD} was proposed and analysed in \cite{EEWZ}.  In that model the cell density cannot exceed a pregiven threshold, extreme aggregation   is avoided, and the propagation speed is finite.
However, as observed in \cref{SecDiff} for the purely diffusion case,  a density-independent upper bound for the  propagation speed and the reproduction of  experimentally observed sharp moving fronts  cannot be achieved in this manner. Similar to the pure diffusion case (see \cref{SecDiff}), this motivates the use of  FL mechanisms.

 FL taxis models rely on replacing $\chi c\nabla_x S$ with 
 \begin{align*}
    J_{chemo}=cV_{chemo},\qquad V_{chemo}=b_{chemo}(c,S,\nabla_x S),
 \end{align*}
where $b_{chemo}$ is a (vector-valued) function such that on the one hand,  $b_{chemo}(c,S,y)$ is close to $\chi(c,S)y$ for  sufficiently small $|y|$  for some bounded function $\chi$, but on the other hand,  $b_{chemo}(c,S,\cdot)$ is bounded for every fixed pair $(c,S)$. The latter property ensures a limitation of the taxis component of the flux.  A prototypical model with FL taxis is thus
\begin{subequations}\label{KSFS}
 \begin{align}
  &\partial_t c=\nabla_x\cdot\left(D_c\nabla_x c-cb_{chemo}(c,S,\nabla_x S)\right),\label{KSFSc}\\
  &\partial_t S=D_v\Delta_x S- \alpha S+ \beta c.\label{KSFSv}
 \end{align}
\end{subequations}
%
 One could, for instance, use the following function that was proposed in \cite{HillenPainter}:
\begin{align*}
 b_{chemo}(c,S,y)=\chi C\left(\tanh\left(\frac{y_1}{1+C}\right),\dots,\tanh\left(\frac{y_d}{1+C}\right)\right)
\end{align*}
for some constants $C$ and $\chi$. 
 We refer to \cite{HillenPainter} and references therein as well as to  \cite{CKWW,PerthameVW,BBNS2010} for further examples of  models with non-FL diffusion and FL chemotaxis. In \cite{KimFried} a model for glioma invasion was developed which includes  linear diffusion and FL chemo- and haptotaxis, as well as other relevant effects.
 

In \cite{BBNS2010} a model with a fully limited cell flux was proposed:
\begin{subequations}\label{new}
 \begin{align}
  &\partial_t c=\nabla_x\cdot\left(D_c\frac{c\nabla_x c}{\sqrt{c^2+\frac{D_c^2}{C^2}|\nabla_x c|^2}}-\chi c \frac{\nabla_x S}{\sqrt{1+|\nabla_x S|^2}}\right)+f_c(c,S),\label{newc}\\
  &\partial_t S=D_v\Delta S+f_v(c,S),\label{newv}
 \end{align}
 \end{subequations} 
 with some functions $f_c$ and $f_v$ and $C>0$ a constant. Both motion effects in \cref{newc} reflect some sort of optimal transport \cite{BBNS2010}. The diffusion term originates from the relativistic heat equation \cref{RH}. Most importantly, choosing both  diffusion and taxis FL guaranties that the speed of propagation cannot exceed a universal constant, in this case $C+\chi$.

 Travelling wave analysis for some  parabolic-elliptic systems with FL diffusion and non-FL chemotaxis \cite{arias2018cross,CYPerthame,campos2021kinks}  and numerical simulations for an extension of \cref{new} to a model for glioblastoma (the most aggressive type of glioma)  invasion that includes FL diffusion as well as multiple FL taxis terms \cite{conte2021modeling} indicate the ability of such  models to propagate singularities observed in biological applications, including cancer invasion. Further models for glioma invasion which involve  FL motion terms  were developed in  \cite{KumarSur,DKSS20}. 
 \bigskip
 
 Similar to the purely diffusion case,  RDT systems with FL mechanisms can be constructed directly on the macroscale. More accurate derivations based on a multiscale approach are discussed in \cref{Sec:Multi} below.

\subsection{Analytical challenge} 
FL diffusion and taxis terms may have a similar form, as e.g. in  \cref{new}, yet their impact on the analysis is vastly different. 
In order to see that, let us compare  \cref{KSFS} and \cref{new} with the classical model \cref{KS}. We assume that these three systems are stated in a bounded domain and we impose the no-flux boundary conditions. 

Classical theory \cite{Amann1}  implies that \cref{KS} is uniquely solvable in the classical sense as long as it remains bounded. If a solution blows up at some finite time, then it still exists globally in a certain generalised sense \cite{ZhigunKS}. Key to solvability is in both cases the special structure of \cref{KS}: it belongs to the class of regular quasilinear upper-triangular parabolic systems  \cite{Amann1}. In such systems, diffusion is sufficiently strong compared to taxis. This allows to obtain certain necessary a priori estimates by manipulating both equations. Thanks to these estimates, existence and uniqueness of maximal classical solutions can be obtained by means of a standard argument which is based on the Banach fixed-point theorem.

In general, systems with linear diffusion and FL taxis can be handled very similar to \cref{KS}, see e.g. \cite{PerthameVW,CKWW}. In  \cref{KSFS} flux-limitation ensures that cell diffusion is the dominating factor in \cref{KSFSc}, which considerably simplifies the analysis.  For example, since the velocity component due to taxis is bounded by construction, a priori boundedness of  $c$ can be obtained by dealing with equation \cref{KSFSc} alone.

On the other hand, system \cref{new} raises new  challenges compared to  \cref{KS}. The FL diffusion term precludes the application of the standard theory of parabolic PDEs even if global a priori boundedness is guaranteed. For a single purely diffusion equation such as \cref{RH} it was possible to establish the existence of a mild  solution and to prove that this solution is also the unique entropy solution \cite{CCCSS15}. For the strongly coupled chemotaxis system \cref{new} it seems that it is neither possible to  set up a semigroup, nor to prove uniqueness even if entropy inequalities are imposed.
A rigorous analysis of 
this system 
is still lacking. 

Currently available analytical studies of systems with a fully FL cell flux are restricted to parabolic-elliptic versions. They are generally easier to deal with than the parabolic-parabolic ones.
In    \cite{BelWinkler20171,BelWinkler20172,MIZUKAMI2019,Chiyoda2019} existence and blow-up were addressed in the radial-symmetric case and for strictly positive initial $c$-values. Travelling wave analysis in \cite{arias2018cross,CYPerthame,campos2021kinks} allows for biologically relevant nonnegative  densities.

\section{Multiscale derivations}\label{Sec:Multi}
In this Section we turn to  derivations which start with a  KTE on the mesoscale and yield a   macroscopic PDE that contains FL diffusion and/or taxis.  This multiscale approach makes possible a more careful modelling than a single-scale purely macroscopic one. 

RDT systems, such as those discussed in \cref{SecTax},  are among  the most widely used tools in cancer invasion modelling. They describe the evolution of macroscopic densities of cancer cell populations and densities/concentrations of other involved components, such as, e.g. tissue or  biochemical signals. These quantities depend only upon time and position in space, which  allows for comparison with  information acquired by  standard biomedical imaging techniques,  e.g. magnetic resonance imaging (MRI)  and computed tomography (CT). Another advantage is the availability of well-developed mathematical analysis tools  and efficient numerical methods.

Often, macroscopic RDT systems are derived using a standard single-scale approach based on the balancing of macroscopic fluxes. Derivations of this type which focus on cell-tissue interaction in cancer  include, e.g.  \cite{anderson2000,ChapLol2011}. However, when modelling directly on the microscale, one may lose important  lower-level information or capture it inaccurately. In contrast, a multiscale modelling approach begins with putting together equations for processes of which at least some are occurring on  scales smaller than the macroscale, i.e. micro, meso, etc. A more or less realistic setting typically includes a  combination of several scales. Unfortunately,  the resulting detailed equations are generally too difficult to solve numerically.
Hence, a suitable upscaling is usually performed,  yielding an RDT system which  still contains some essential lower-level information in the equation terms. If successfully studied analytically and simulated numerically, using parameters determined from experiments, it generally offers a much more careful description of cell migration than that which can be achieved through   balancing of  (macroscopic) fluxes.

The described multiscale modelling approach was originally applied in physical context. It allowed, for example, to obtain the Euler and  Navier-Stokes equations as macroscopic scaling limits of the Boltzmann equation, see e.g. \cite{SaintRaymond}. 
Later, the approach based on modelling with KTEs and their subsequent upscaling was successfully adjusted to modelling of population motility. Starting from \cite{Alt1980,OthmerDunbarAlt} numerous models have been derived in this manner. 

In the context of cancer migration, the KTE-based approach  allows for an adequate description of the impact that various sorts of  heterogeneity have on tumour invasion. For instance, the cell-tissue interaction depends on such variables as: the tissue fiber position-direction  distribution (mesoscopic), the amount of free receptors on the cell surface  that can bind to tissue (microscopic), the density/concentration gradients of various tactic  signals  (macroscopic), etc.
Including such variables leads to multiscale settings which, when upscaled, result in nonstandard PDEs that  differ considerably from those set directly on the macroscale. We refer, e.g. to \cite{HillenPainter} (see also references therein), where the effect of environmental anisotropy is comprehensively addressed (not specifically for cancer), showing that it leads to equations with drift  and/or myopic (and thus non-Fickian) diffusion, both of which depend on  parameters of the tissue distribution. The kinetic theory for active particles (KTAP) \cite{bellom3} is a further development of  the method for those situations where not only physical variables (time, position, velocity, etc.) but, also, the so-called 'active variables' are involved. For example,  in \cite{KelkelSur2012} an  extension of an earlier model for glioma invasion  under tissue anisotropy \cite{HillenM5}  is presented which treats cell surface receptors as active variables. That model also includes chemo- and haptotaxis.

\bigskip

To simplify the exposition, we start with a single KTE of the form
\begin{align}
\partial_t c+v\cdot\nabla_x c={\mathcal L}(c)\label{KTE}
\end{align}
for  cell density distribution $c$. This  mesoscopic quantity is a function not only of time $t\geq0$ and position $x\in\R^d$, $d\in\N$ but also of  velocity $v$ which belongs to a suitably chosen bounded
velocity space ${\mathcal V}\subset \R^d$. Equation \cref{KTE} balances  deterministic transport and probabilistic changes.  For a completely zero right-hand side, we obtain a simple transport equation  corresponding to the situation when  velocities of cells do not change. On the microscopic level, the movement of each of them is modelled by the ODE system
\begin{subequations}\label{MicroDSt}
\begin{align}
&\frac{dx}{dt}=v, \\
&\frac{dv}{dt}=0,\label{eq:microdynSt}
\end{align}
\end{subequations}
Of course, this is only valid if no deterministic external forces are acting on the cells. Later in \cref{withExt} we  consider a case where such forces are present.  We assume the operator ${\mathcal L}$ on the right-hand side of \cref{KTE} to be a so-called turning operator: it models the impact of a velocity-jump process, i.e. of  probabilistic instantaneous changes in the velocity of the species. 

Aiming at a macroscopic description, one often performs a suitable upscaling. The standard approach  begins with a rescaled equation
\begin{align}
 \varepsilon^{\kappa} \partial_t \ce+\varepsilon v\cdot\nabla_x \ce={\mathcal L}_{\ve}(\ce),\label{KTEe}
\end{align}
where $t$ and $x$ now stand for macroscopic time and space variables, respectively,  $\varepsilon>0$ is a small scaling parameter, and, typically, $\kappa=2$ or $\kappa=1$, corresponding to the parabolic and hyperbolic scalings, respectively.  In order to obtain a macroscopic characterisation, one seeks to eliminate $v$ and looks for a good approximation of $\overline \ce=\overline \ce(t,x)$ as $\ve$ tends to zero. Here and in what follows we use the following notation.
\begin{Notation}
We denote 
$$\overline u\overset{\text{def}}{=}\int_{{\mathcal V}} u\,dv$$ 
if $u:{\mathcal V}\to\R$ is a function and
$$\overline u\overset{\text{def}}{=}\int_{{\mathcal V}} du$$
if it is a measure in ${\mathcal V}$. 
\end{Notation}
As a rule, one assumes that $\ce$ can be well-approximated by a truncation of the Hilbert expansion
\begin{align}
 \ce=\sum_{n=0}^{\infty}\ve^n c_n^0,\label{Hilbert} 
\end{align}
where $\overline{c_0^0}$ is of particular importance since it is the zero-order macroscopic approximation. The first-order correction $\overline{c_1^0}$ is also of interest.

Since  integration and a nonlinear map do not commute, it is in general difficult, if not impossible, to deal with  a nonlinear ${\mathcal L}$ unless it is, e.g. of the form 
\begin{align}
{\mathcal L}(c)={\mathcal L}[f_1,\dots,f_M]c,\end{align}
 for some macroscopic  functions $f_i=f_i(t,x)$, $i=1,\dots,M$, $M\in\N$,  which when fixed define a linear operator  ${\mathcal L}[f_1,\dots,f_M]$. It is possible to have $\overline c$ among $f_i$'s. 
\bigskip

Three approaches  starting  from such equations as \cref{KTE} or its extensions and leading to RDT systems with FL effects have been proposed  so far.  We briefly review them in the reminder of this Section.

\subsection{Nonlinear Hilbert expansion}\label{NLHil}
In this Subsection we look at a construction where FL diffusion results from a nonstandard approximation.
In \cite{Coulombel2005}, one considered a parabolic scaling of \cref{KTE} for the turning operator which, in the above notation, takes  the form
\begin{align}
 {\mathcal L}c=\lambda(\overline{c}\mu-c),\label{OperL}
\end{align}
where $\lambda>0$ is a constant, $\mu$ is a fixed probability measure  on a bounded space ${\mathcal V}\subset[-1,1]$ such that for any continuous odd function $h$
\begin{align}
 \int_{{\mathcal V}}h\,d\mu=0\label{vmu}
\end{align}
and satisfying certain other  assumptions, and $c=f\mu$ for a square integrable density $f$. 
As set out in \cite{HillenPainter}, this kind of turning operator is a convenient way to capture gains and losses due to instantaneous velocity changes for a population moving in a heterogeneous surroundings. In the context of cancer migration, $\mu$  typically stands for the orientational distribution of tissue
fibers, see e.g. \cite{HillenM5,KelkelSur2012}, whereas $\lambda$ is the reorientation rate. 

Since ${\mathcal L}$ is linear, the corresponding equations for $\overline{c_0^0}$ and $\overline{c_0^0}+\ve\overline{c_1^0}$ are linear as well. In fact, both approximations satisfy the linear diffusion equation \cref{LD} with $$D_c=\frac{1}{\lambda}\int_{{\mathcal V}}v^2\,d\mu,$$
so that the propagation speed is infinite (see the discussion in \cref{SecDiff}). Yet this is not the case for $\ce$. Indeed, integrating \cref{KTEe} yields the  conservation law
\begin{align*}
 \partial_t \overline{\ce}+\partial_x(\overline{\ce} V^{\ve})=0, \qquad V^{\ve}=\frac{1}{\ve}\frac{\int_{{\mathcal V}}v\,d{\ce}}{\int_{{\mathcal V}}d{\ce}},
\end{align*}
and, since $V\subset [-1,1]$,
\begin{align*}
 |V^{\ve}|\leq \frac{1}{\ve}<\infty.
\end{align*}
To avoid  the infinite propagation speed on the macroscopic scale, as well as further undesirable effects, it was proposed in  \cite{Coulombel2005} to consider a nonlinear Hilbert expansion, replacing  \cref{Hilbert} with 
\begin{align}
 \ce=e^{\sum_{n=0}^{\infty}\ve^n \Phi_n^0}\mu.\label{Hilbertexp} 
\end{align}
It was proved there that a  good approximation of the  first order truncation 
\begin{align}
 \overline{e^{\Phi_0^0+\ve \Phi_1^0}\mu}\nonumber
\end{align}
and, thus, of $\ce$  solves
\begin{align}
 \partial_t u^{\ve}=\partial_x\left(\frac{1}{\ve}u^{\ve}\G\left(\frac{\ve}{\lambda}\frac{\partial_xu^{\ve}}{u^{\ve}}\right)\right),\label{FSG}
\end{align}
where 
\begin{align}
 \G(\beta)=&\frac{\int_{{\mathcal V}}ve^{\beta v}\,d\mu}{\int_{{\mathcal V}}e^{\beta v}\,d\mu}\nonumber\\
 =&\frac{d}{d\beta}\ln \left(\int_{{\mathcal V}}e^{\beta v}\,d\mu\right),\label{G}
\end{align}
and, under the assumptions on $V$ and $\mu$ as imposed in  \cite{Coulombel2005},
\begin{enumerate}
 \item\label{G1} $\G:\R\to(-1,1)$ is odd, strictly increasing, infinitely differentiable, and a diffeomorphism;
 \item\label{G2} $\G(\pm\infty)=\pm1$.
\end{enumerate}
In particular, the saturation property is guaranteed. 
For some measures the corresponding $\G$ can be computed explicitly. 
\begin{Example}[\cite{Coulombel2005}]\label{Examplemu}
Let ${\mathcal V}=[-1,1]$.
 \begin{itemize}
 \item A homogeneous environment corresponds to the normalised Lebesgue measure:
\begin{align}
 &\mu=\frac{1}{2}|\cdot|\nonumber
 \\
\Rightarrow\quad&\G(\beta)=\coth(\beta)-\frac{1}{\beta}.\nonumber
\end{align}    
\item If cell speeds remain constant, then one deals with a discrete measure:
\begin{align}
 &\mu=\frac{1}{2}(\delta_{-1}+\delta_1)\label{mud}\\
\Rightarrow\quad&\G(\beta)=\tanh(\beta).\nonumber
\end{align}    
\end{itemize}

\end{Example}
However, not every $\G$ which satisfies \cref{G1}-\cref{G2} can be generated in this way.  
To see why this is indeed the case, let us consider
\begin{align*}
 I(\beta)\overset{\text{def}}{=}e^{\int_0^{\beta}G(s)\,ds}= \int_{{\mathcal V}}e^{\beta v}\,d\mu,
\end{align*}
where the latter equality is due to \cref{G}. Then, for all $n\in\N$
\begin{align}
 \frac{d^{2n}}{d\beta^{2n}}I(0)= \int_{{\mathcal V}}v^{2n}\,d\mu\in(0,1].\label{an}
\end{align}
Many interesting $\G$'s do not satisfy \cref{an}. This includes a function which corresponds to the relativistic heat equation: 
$$\G(\beta)=\frac{\beta}{\sqrt{1+|\beta|^2}}.$$ Indeed, a direct computation shows that
\begin{align*}
 I(\beta)=e^{\int_0^{\beta}G(s)\,ds}=e^{\sqrt{1+|\beta|^2}-1}
\end{align*}
and 
\begin{align*}
 \frac{d^{4}}{d\beta^{4}}I(0)=0,
\end{align*}
which violates \cref{an}.
\bigskip

Yet another limitation is that in general the resulting FL diffusion equation \cref{FSG} generates the propagation speed  $1/\ve$ (see \cref{SecDiff}) which, while finite, may be unrealistically  large.
\subsubsection*{Summary}
\begin{enumerate}
 \item A parabolic scaling of a  basic linear KTE together with a  nonlinear 'exponential' Hilbert expansion leads to a FL (and thus nonlinear) diffusion equation.
 \item This construction cannot be used in order to obtain some standard FL equations, including the relativistic heat equation.
 \item The resulting generic propagation speed is $1/\ve$.
\end{enumerate}

\subsection{Scaled  turning operator}\label{Sec:scaledOper}
We have seen in the preceding Subsection that a FL effect (diffusion in that particular  case) can be obtained from a basic linear KTE if one is prepared to go beyond the zero-order approximation, while keeping some sort of first order correction. This is often  undesirable, i.e. one wishes to have a macroscopic  PDE for $\overline{c_0^0}$ alone. 
Studies in  \cite{BBNS2010,PerthameVW,DolakSchmeiser} show that it is possible to obtain  equations featuring  FL diffusion and/or taxis by  considering operators of   \cref{OperL} type and choosing  $f_i$'s and their scalings appropriately. We exploit such constructions in this Subsection. 

To illustrate the approach, we use a turning operator of the form
\begin{align}
&{\mathcal L}[S]c={\mathcal L}_0c+{\mathcal L}_1[S]c,\nonumber\\
&{\mathcal L}_0c=\lambda\left(\frac{1}{|{\mathcal V}|}\overline c-c\right),\nonumber\\
&{\mathcal L}_1[S]c=\int_{\mathcal V}T_1[S](v,v')c(v')-T_1[S](v',v)c(v)\,dv'.
\end{align}
It models a superposition of two essentially independent effects. 
For $\lambda>0$, the first component, ${\mathcal L}_0$ is the multidimensional version of a special case of \cref{OperL} for the normalised Lebesgue measure $$\mu=\frac{1}{|{\mathcal V}|}|\cdot|.$$  
This choice corresponds to  chaotic velocity changes in a homogeneous environment and amounts to a constant linear diffusion operator in the equation for the macroscopic zero order approximation $\overline{c^0}$.  The second component of ${\mathcal L}[S]$, operator ${\mathcal L}_1[S]$  has been added with the aim to capture  velocity changes due to the influence of a substance with density/concentration $S=S(t,x)$.  For each velocity pair $v',v\in {\mathcal V}$, the corresponding value $T_1[S](v,v')$ of a so-called turning kernel $T_1[S]$ can be interpreted as the likelihood of a cell to change from $v'$ to $v$, provided that $T_1[S]$ is  nonnegative. The kernel needs to satisfy certain further assumptions, so that, in particular, ${\mathcal L}_1[S]$ is a conservative operator, i.e. 
\begin{align*}
 \int_{{\mathcal V}}{\mathcal L}_1[S]c\,dv=0.
\end{align*}
Such turning operators are a standard tool for deriving chemotaxis models, with $S$ being the concentration of a chemoattractant, see e.g.  \cite{BBTW} and references therein.   Here we concentrate solely on equations for cell dynamics since the dynamics of the chemical is a standard macroscopic one. 

In the reminder of this subsection, we review the effect of some possible choices of $T_1[S]$ and of their scalings such that lead to RDT equations with FL diffusion and/or taxis terms on the macroscale. 
\subsubsection{Dependence on  past motion.} \label{SecPast}
We begin with 
\begin{align}
 T_1[S](v,v')=\Psi(D_tS),\qquad D_tS=\partial_tS+v'\cdot\nabla_x S.\label{T1}
\end{align}
Proposed in \cite{DolakSchmeiser}, it describes the likelihood of a velocity change as a function of the temporal derivative of $S$ along the path the cell has been moving prior to that change. This choice is based on the known ability of cells to compare present  signal concentrations  to previous ones and respond to that.  Function $\Psi$ describes the response rate. We assume it to be nonlinear.

Different upscalings can be adopted for \cref{KTE} with the turning kernel \cref{T1}. The hyperbolic limit is a drift equation for $\overline{c^0}$ \cite{DolakSchmeiser}. No diffusion can be recovered this way unless a first-order correction is included. 
A straightforward parabolic rescaling would be
\begin{align}
 &\ve^2\partial_t \ce+\ve v\cdot\nabla_x \ce\nonumber\\
 =&\lambda\left(\frac{1}{|{\mathcal V}|}\overline \ce-\ce\right)\nonumber\\
 &+\int_{\mathcal V}\Psi\left(\ve^2\partial_tS^{\ve}+\ve v'\cdot\nabla_x S^{\ve}\right)\ce(v')-\Psi\left(\ve^2\partial_tS^{\ve}+\ve v\cdot\nabla_x S^{\ve}\right)\ce(v)\,dv'.\nonumber
\end{align}
However, a nonlinear dependence upon the gradient of the attractant is lost in the limit when $\ve$ is sent to zero. To preclude this, the response function needs to be rescaled as well. In \cite{PerthameVW}, one therefore replaced $\Psi$ by
$$\ve\Psi\left(\frac{\cdot}{\ve}\right).$$
This choice is biologically justifiable, see \cite{PerthameVW} and references therein. The resulting  rescaled KTE is 
\begin{align}
& \ve^2\partial_t \ce+\ve v\cdot\nabla_x \ce\nonumber\\
 =&\lambda\left(\frac{1}{|{\mathcal V}|}\overline \ce-\ce\right)\nonumber\\
 &+\ve\int_{\mathcal V}\Psi\left(\ve\partial_tS^{\ve}+ v'\cdot\nabla_x S^{\ve}\right)\ce(v')-\Psi\left(\ve\partial_tS^{\ve}+v\cdot\nabla_x S^{\ve}\right)\ce(v)\,dv'.\label{KTETaxe1}
\end{align}
It was verified in \cite{PerthameVW} that this leads to the diffusion-taxis equation
\begin{align}
 \partial_t\overline{c^0}=D\Delta_x\overline{c^0}-\nabla\cdot(\overline{c^0}\Phi(\nabla_x S^0))\label{MacroFSC}
\end{align}
where
\begin{align}
 &D=\frac{1}{\lambda|{\mathcal V}|}\int_{{\mathcal V}}v\otimes v\,dv,\nonumber\\
 & \Phi(\beta)=-\frac{1}{\lambda}\int_{{\mathcal V}}v\Psi(v\cdot \beta)\,dv,\label{PhiPsi}
\end{align}
and, under suitable conditions on $\Psi$, $\ce$ is well-approximated by $c^0$.  

Various kinds of chemotactic response can be obtained in this manner.
Let us consider a basic case: $${\mathcal V}=[-1,1].$$ Relation  \cref{PhiPsi} between $\Phi$ and $\Psi$ implies that  $\Phi$ is necessarily odd and that then  
\begin{align*}
 \Psi(\beta)=\Psi_{even}(\beta)-\frac{\lambda}{2\beta}\frac{d}{d\beta}(\beta^2\Phi(\beta)),
\end{align*}
where $\Psi_{even}$ is any even function.
In particular, a direct computation shows that choosing 
\begin{align}
 \Psi(\beta)=C-\frac{x \left(2 x^2+3\right)}{2\left(x^2+1\right)^{3/2}},\label{PsiEx}
\end{align}
with $C$ a constant, leads to 
\begin{align*}
\Phi(\beta)=\frac{\beta}{\sqrt{1+|\beta|^2}},
\end{align*}
 so that the corresponding taxis term is as in \cref{new}. Choosing $C$ sufficiently large ensures that $\Psi$ is nonnegative and can therefore be viewed as the likelihood of turning due to taxis. However, it is not clear how to interpret the particular form of $\Psi$ in  \cref{PsiEx}.  Further examples of possible $\Psi$ can be found, e.g. in \cite{DolakSchmeiser} and  \cite{PerthameVW} (see also references therein). 

The approach works for various modifications of the turning operator ${\mathcal L}$. For instance, the same equation \cref{MacroFSC} is obtained if one uses
\begin{align*}
 &T_1[S](v,v')=\Psi(v'\cdot\nabla_x S)
\end{align*}
instead of \cref{T1} since $\partial_tS$ only appears on the macroscale if a first order correction is included. 
In both  cases one could take $S$ to be a function of $\overline{c}$, such as e.g.
$$S=-\ln(\overline c).$$ This  leads to a  diffusion flux which is nonlinear with respect to $\nabla_x\overline{c^0}$, yet it is not FS due to the Fickian  contribution from ${\mathcal L}_0$.  

One could also use more general ${\mathcal L}_0$, e.g. such as in \cref{OperL} for a non-Lebesgue measure in order to account for the environmental heterogeneity. 
 \subsubsection{Dependence on both anterior and posterior velocities.}\label{Sec:FullFS}
 Aiming at a fully FL cell flux, 
  the following kernel was proposed in \cite{BBNS2010}:
 \begin{align}
  &T_1[\alpha[\overline c,S]](v,v')=\left(\alpha[\overline c,S]-v'\right)\cdot vh(v),\label{T1al}\\
  &\alpha[\overline c,S]=D_c\frac{\nabla_x \overline c}{\sqrt{\overline c^2+\frac{D_c^2}{C^2}|\nabla_x \overline c|^2}}-\chi  \frac{\nabla_x S}{\sqrt{1+|\nabla_x S|^2}},\nonumber
 \end{align}
 where  $h$  satisfies 
 \begin{align}
  \int_{{\mathcal V}}h(v)\,dv=1,\qquad \int_{{\mathcal V}}vh(v)\,dv=0,   \qquad \int_{{\mathcal V}}v\otimes vh(v)\,dv=\beta I \nonumber
 \end{align}
  for a positive constant $\beta$, and $S$ is the concentration of a signal substance.  
The authors then considered an extension of \cref{KTE} involving further  integral operators, e.g.  such that model cell proliferation on the mesoscale. They took 
\begin{align*}
 \lambda=0
\end{align*}
and considered the KTE
\begin{align}
 \varepsilon \partial_t \ce+\varepsilon v\cdot\nabla_x \ce={\mathcal L}[\alpha[\overline{\ce}],S^{\ve}]]\ce+\text{[growth, etc.]}.\label{KTEeh}
\end{align}
This equation can be interpreted as a hyperbolic scaling for \cref{KTE} with a rescaled turning kernel \cref{T1al}: $\alpha[\overline c,S]$ needs to be replaced by
 \begin{align*}
  D_c\frac{\frac{1}{\ve}\nabla_x \overline c}{\sqrt{\overline c^2+\frac{D_c^2}{C^2}\left|\frac{1}{\ve}\nabla_x \overline c\right|^2}}-\chi  \frac{\frac{1}{\ve}\nabla_x S}{\sqrt{1+\left|\frac{1}{\ve}\nabla_x S\right|^2}}
 \end{align*}
in order to have  \cref{KTEeh} after a hyperbolic scaling.
Proceeding with a formal limit as $\ve\to 0$, one recovered in \cite{BBNS2010} the fully FL diffusion-taxis equation \cref{newc}.
\bigskip

While choosing the kernel as in \cref{T1al} ensured the desired macroscopic limit, one has that

\begin{enumerate}
\item  the corresponding integral operator ${\mathcal L}_1[\alpha[\overline c,S]]$ is conservative,
\item but the kernel is not nonnegative everywhere.                                                                        \end{enumerate}
Thus, ${\mathcal L}_1[\alpha[\overline c,S]]$ is not a turning operator in the traditional sense. 

\subsubsection*{Summary}
 \begin{enumerate}
\item 
It is possible to obtain an RDT equation with FL effects, including a completely FL diffusion-taxis flux, as the zero order approximation of KTE \cref{KTE}  by choosing suitably:
\begin{itemize}
\item[(i)] a turning kernel,
\item[(ii)] a scaling of time and space (parabolic/hyperbolic),
\item[(iii)] a rescaling of the turning kernel.
\end{itemize}      
\item The taxis term on the macroscale arises from a turning operator and thus has probabilistic roots.
\item The approach is flexible, but may require using a turning kernel that is difficult to interpret, e.g. because it has a complicated form and/or it is not everywhere nonnegative.
\end{enumerate}

\subsection{Accelerated motion}\label{withExt}
So far in this Section we have dealt with constructions which presuppose that cell velocity changes are fully probabilistic.  In particular, the FL effects on the macroscale stemmed from the turning kernel. 
In this final Subsection we turn to a different approach that was developed in \cite{DSW,ZSMM}. As an illustration we use a special case  of the modelling framework from \cite{ZSMM}. 
 There the deterministic part of cell motion is described on the microscale by the following extension of \cref{MicroDSt}:
\begin{subequations}\label{MicroD}
\begin{align}
&\frac{dx}{dt}=v, \\
&\frac{dv}{dt}=-a(v-v_*[S](t,x)),\label{eq:microdyn}
\end{align}
\end{subequations}
with
\begin{align}
v_*[S]=\F\frac{\nabla_x S}{1+|\nabla_x S|}.\label{vstar}
\end{align}
The ODE system \cref{MicroD} resembles the second Newton's law. However, unlike lifeless matter for which the acceleration would necessarily be due to a physical force, here  the presence of a signal  stimulates the cells to divert from a straight line. 
The choice of the right-hand side in \cref{eq:microdyn} is motivated by the assumption that a cell tends to realign with a certain  'preferred' velocity $v_*[S]$ which depends on the spacial gradient of a cue with density/concentration $S=S(t,x)$  and  on the spacial 
heterogeneity of the environment. The latter is accounted for by means of  a matrix-valued function $\F=\F(x)$. In the absence of signal gradients cell deceleration is proportional to a constant $a>0$, as in the Stokes' law. 
Choosing $$\|\F(x)\|_2\leq 1\qquad\text{for all }x\in\R$$ ensures that $v_*[S]$ remains inside the bounded velocity space 
\begin{align*}
 {\mathcal V}=\{x\in\R^d:\ \ |x|<1\}.
 \end{align*}
 Combining the microscopic dynamic \cref{MicroD}  with a turning operator which, for simplicity, we assume here to be of the form
 \begin{align*}
  {\mathcal L}c=\lambda\left(\frac{1}{|{\mathcal V}|}\overline c-c\right)
 \end{align*}
for some constant $\lambda>0$, and the mass conservation law, one arrives at the KTE
\begin{align}
 \partial_t c+\nabla_x\cdot (vc)-a\nabla_v\cdot((v-v_*)c)
 =&\overline{c}-c.\label{meso}
\end{align}
The parabolic scaling then yields a diffusion-taxis equation for the zero order approximation:
\begin{align}
 (a+\lambda)\partial_t \overline{c^0}
 =\frac{\lambda}{2a+\lambda}\frac{n}{n+2}\Delta_x\overline{c^0}-a\nabla_x\cdot\left(\overline{c^0}\F\nabla_x S\right),\label{dif2}
\end{align}
and it was rigorously proved in \cite{ZSMM} that $\ce$ is well-approximated by   $c^0$. While in  derivations in  \cref{Sec:scaledOper} all  motions effects on the macroscale came from a turning operator, here it is a source of diffusion only. This time, taxis originates from the deterministic microscopic dynamics. A possible interpretation is that cells change their velocities in the attempt to follow the attractant gradients but may at the same time be diverted from such preferred trajectories by chaotic velocity jumps. 

Alike the derivation in \cref{NLHil}, including a first order correction leads to a FL effect on the macroscale \cite{ZSMM}:
\begin{align}
 (a+\lambda)\partial_t \overline{c_{01}^{\ve}}
 =\frac{\lambda}{2a+\lambda}\frac{n}{n+2}\Delta_x\overline{c_{01}^{\ve}}-a\nabla_x\cdot\left(\overline{c_{01}^{\ve}}\F\frac{\nabla_x S}{1+\ve|\nabla_x S|}\right)+O\left(\ve^2\right)\label{dif3}
\end{align}
for the first order approximation
\begin{align*}
 c_{01}^{\ve}=c^{0}+\ve c^{0}_1.
\end{align*}
If the error term on the right-hand side is dropped, \cref{dif3} becomes a parabolic PDE with a linear diffusion and a FL taxis. However, and this is also related to what we saw in \cref{NLHil}, the velocity component due to taxis is of order $O(1/\ve)$, i.e. potentially too large. 
Similar to \cref{Sec:FullFS}, one could avoid $\ve$ in the limit equation by replacing $v_*[S]$ in \cref{eq:microdyn} with 
\begin{align}
v_*^{\ve}[S]=\F\frac{\frac{1}{\ve}\nabla_x S}{1+\frac{1}{\ve}|\nabla_x S|}\nonumber
\end{align}
and adopting the hyperbolic scaling of the KTE.
One readily verifies that this amounts to 
\begin{align}
 (a+\lambda)\partial_t \overline{c^0}=a\nabla_x\cdot(-\overline{c^0}v_*[S^0]).\label{hyp}
\end{align}
A similar scaling was performed in \cite{DSW} in a more general context. There one first used the moments method and then the hyperbolic scaling. 
In \cite{KumarSur} one relied fully on the moment closure, so no rescaling was required.

Different forms of $v_*[S]$ and more general turning operators are  possible. For example, for 
\begin{align*}
&v_*[\overline c]=\psi\left(-\nabla_x \ln(\overline c)\right),\\
 &v_*^{\ve}[\overline c]=\psi\left(-\frac{1}{\ve}\nabla_x \ln(\overline c)\right)
\end{align*}
with $\psi$ as in \cref{psi} the same procedure leads to  equation \cref{hyp} which now  becomes the relativistic heat equation \cref{RH}. Choosing the turning operator ${\mathcal L}$ as, e.g. in \cref{OperL} for a non-Lebesgue measure leads to an anisotropic diffusion and  additional transport terms. This was done in \cite{ZSMM,DSW,KumarSur} in order to account for a fibrous  environment. While \cite{ZSMM} dealt with scalings and their rigorous justification for a single prototypical equation,  \cite{DSW,KumarSur} focused on applications, providing  detailed models for tumour invasion.  

\subsubsection*{Summary}
\begin{enumerate}
 \item Including a suitable transport term with respect to velocity into the mesoscopic KTE allows to obtain an RDT equation with FL effects.
 \item The taxis term on the macroscale arises from the  deterministic motion component and is independent from the turning kernel.
 \item The approach is flexible, but may require a suitable rescaling of terms in order to avoid
unrealistically large propagation speeds.
 \end{enumerate}
\begin{acknowledgement}
The results of this paper were presented by the author at the 34th Annual Meeting of the Irish Mathematical Society in September 2021. The author was supported by the Engineering and Physical Sciences Research Council [grant number EP/T03131X/1]. For the purpose of open access, the author has applied a Creative Commons Attribution (CC BY) licence to any Author Accepted Manuscript version arising.
\end{acknowledgement}

\addcontentsline{toc}{section}{References}
\printbibliography


\end{document}